\documentclass[12pt,oneside]{amsart}
\usepackage{amsmath}
\usepackage[mathscr]{euscript}
\usepackage{epsfig}
\usepackage{graphics}
\usepackage{amsbsy}
\usepackage{amsfonts}
\usepackage{amssymb}
\usepackage{graphicx}
\usepackage{rotating}
\usepackage{cancel}

 \newcommand{\bd}{\begin{definition}}
 \newcommand{\ed}{\end{definition}}
 \newcommand{\bt}{\begin{theorem}}
 \newcommand{\et}{\end{theorem}}
 \newcommand{\bp}{\begin{proposition}}
 \newcommand{\ep}{\end{proposition}}
 \newcommand{\bl}{\begin{lemma}}
 \newcommand{\el}{\end{lemma}}
 \newcommand{\bpr}{\begin{proof}}
 \newcommand{\epr}{\end{proof}}
 \newcommand{\bc}{\begin{corollary}}
 \newcommand{\ec}{\end{corollary}}
\newcommand{\br}{\begin{remark}}
\newcommand{\er}{\end{remark}}

\newcommand{\lra}{\longrightarrow}

\newcommand{\lmt}{\longmapsto}

\renewcommand{\leq}{\leqslant} 
\renewcommand{\geq}{\geqslant}

\newcommand{\gen}[1]{\langle#1\rangle}

\newcommand{\M}{\mathcal{M}}
\newcommand{\T}{\mathbb{T}}
\newcommand{\F}{\mathbb{F}}

\newcommand{\spn}{\mathrm{span}}

\newcommand{\Ima}{\mathrm{Im}}

\usepackage{amssymb}
\usepackage{mathrsfs}
\newtheorem{theorem}{{\bf Theorem}}[section]
\newtheorem{lemma}[theorem]{{\bf Lemma}}
\newtheorem{corollary}[theorem]{{\bf Corollary}}
\newtheorem{proposition}[theorem]{{\bf Proposition}}

\newtheorem{definition}[theorem]{{\bf Definition}}
\newtheorem{remark}[theorem]{{\bf Remark}}
 \makeatletter
      \def\@setcopyright{}
      \def\serieslogo@{}
      \makeatother
\title[The $2$-nilpotent multiplier of $n$-Lie algebras and its applications]{The $2$-nilpotent multiplier of $n$-Lie algebras \\ and its applications}
\normalsize{
\author{Farshid Saeedi$^*$ and Seyedeh Nafiseh Akbarossadat}
}
\begin{document}
\maketitle

\noindent
\textbf{Abstract.} 
In this paper, we first recall the concept of $c$-nilpotent multiplier and $c$-capability of $n$-Lie algebras and also,  recall the formula for calculating the number of basic commutators in $n$-Lie algebras. Then we give the structure of $2$-nilpotent multiplier of the direct sum of two $n$-Lie algebras. Next, we calculate the dimension of $2$-nilpotent multiplier of every abelian $n$-Lie algebras and Heisenberg $n$-Lie algebras $H(n,m)$. Then we give a dimension of $2$-nilpotent multiplier of any nilpotent $n$-Lie algebras of class $2$ by using the number of basic commutators.  
\ \\

\section{Introduction and Preliminary concepts}

\indent 
Lie polynomials appeared at the end of $19${th} century and the beginning of the $20${th} century in the work of Campbell, Baker, and Hausdorff on exponential mapping in a Lie group, which has lead to the so called the  Campbell--Baker--Hausdorff formula. Around $1930$,  Witt introduced the Lie algebra of Lie polynomials and  showed that the Lie algebra of Lie polynomials is actually a free Lie algebra and that its enveloping algebra is the associative algebra of noncommutative polynomials. He proved what is now called the Poincaré--Birkhoff--Witt theorem and showed how a free Lie algebra is related lower central series of the free group. About at the same time,  Hall \cite{M. Hall,Hall} and Magnus \cite{Magnus}, with their commutator calculus, opened the way to bases of the free Lie algebra.  For more details about a historical review of free Lie algebras, we refer the reader to the reference \cite{[21]} and the references therein.

The concept of basic commutators is defined in groups and Lie algebras, and there is also a way to construct and identify them. Moreover, a formula for calculating their number is obtained. 

In 1962,  Shirshov \cite{Shirshov} gave a method that generalizes Hall’s method \cite{M. Hall} for choosing a basis in a free Lie algebra. 

Basic commutators are of particular importance in calculating the dimensions of different spaces and are therefore highly regarded.  Niroomand and Parvizi \cite{Niroomand-Parvizi-M^2(L)} investigated  some more results about $2$-nilpotent multiplier $\mathcal{M}^{(2)}(L)$ of a finite-dimensional nilpotent Lie algebra $L$ and by using the Witt formula, calculated its dimension. Moreover,  Salemkar, Edalatzadeh, and  Araskhan \cite{Salemkar-Edalatzadeh-Araskhan} introduced the concept of $c$-nilpotent multiplier $\mathcal{M}^{(c)}(L)$ of a finite-dimensional Lie algebra $L$ and obtained some bounds for $\mathcal{M}^{(c)}(L)$ by using the Witt formula and basic commutators. 


This paper is organized into three parts.
In the first part, we recall preliminary concepts, such as $c$-nilpotency and $c$-capability,  as well as useful propositions for proving the main results. We also review the formula for counting the number of basis commutators for finite-dimensional $n$-Lie algebras.  In addition, we continue by stating some basic properties and proving that the center and exterior center of any generalized Heisenberg $n$-Lie algebra are equal. 

In the second section, we first give a formula for calculating the $2$-nilpotent multiplier of the direct sum  for $n$-Lie algebras. Furthermore, we also obtain sequences for  $2$-nilpotent multiplier of $n$-Lie algebras. We calculate the dimension  of $c$-nilpotent multiplier for any finite-dimensional  Abelian $n$-Lie algebra, as well as the dimension  of $2$-nilpotent multiplier   of any  Heisenberg $n$-Lie algebra, generalized Heisenberg $n$-Lie algebra, and any nilpotent $n$-Lie algebra  of class $2$ with derivative of dimension $1$.
We present inequalities for the dimension of $2$-nilpotent multiplier of nilpotent $n$-Lie algebras of class $2$ with derivative dimension greater than $1$.

Finally, in the third section, we will examine the conditions for $2$-capability of $n$-Lie algebras.


In 1985, Filippov \cite{vtf} introduced the concept of \textit{$n$-Lie algebras}, as an $n$-ary multilinear  and skew-symmetric operation $[x_1,\ldots,x_n]$ that satisfies the following generalized Jacobi  identity:
\[[[x_1,\ldots,x_n],y_2,\ldots,y_n]=\sum_{i=1}^n[x_1,\ldots,[x_i,y_2,\ldots,y_n],\ldots,x_n].\]
Clearly, such an algebra becomes an ordinary Lie algebra when $n=2$.

Let $L_1,L_2,\ldots,L_n$ be subalgebras of an $n$-Lie algebra $L$. Denote by $[L_1, L_2, \ldots, L_n]$ the subalgebra of $L$ generated by all vectors $[x_1, x_2, \ldots, x_n]$, where $x_i\in L_i$, $i=1, 2, \ldots, n.$ The subalgebra $[L,L, \ldots, L]$ is called the \textit{derived algebra} of $L$, and it is denoted by $L^2$. If $L^2=0$, then $L$ is called an abelian algebra. An \textit{ideal} $I$ of an $n$-Lie algebra $L$ is a subspace of $L$ such that $[I,L, \ldots,L] \subseteq I.$ If $[I, I,L, \ldots ,L]=0$, then $I$ is called an \textit{abelian ideal}.

Let $L$ be an $n$-Lie algebra and let $z\in L$ such that $[z,L,\ldots,L] = 0$. We call the  collection of all such elements in $L$ the center of $L$ and denote it by $Z(L)$. One can check that $Z(L)$ is an ideal of $L$. Put $Z_0(L)=Z(L)$, and define $Z(L/Z_{i-1}(L))=Z_i(L)/Z_{i-1}(L)$. Hence we have $Z_{i-1}(L)\unlhd Z_i(L)$, for all $i\in\Bbb N$. So we can make the following chain: 
\[0\unlhd Z_0(L)\unlhd Z_1(L)\unlhd Z_2(L)\unlhd\cdots\unlhd Z_i(L)\unlhd\cdots.\]
The above series is known as the upper central series.  

Another well-known chain is the lower central series; Let $L^1=L$, let $L^2=[L,L,\ldots,L,L]$, and let $L^{i+1}=[L^{i},L,\ldots,L,L]$. Then we have $L^{i+1}\unlhd L^i$, for $i\geq2$. Thus the following series, known as the lower central series, can be formed as 
\[\cdots\subseteq L^i\subseteq L^{i-1}\subseteq \cdots\subseteq L^3\subseteq L^2\subseteq L^1=L.\]
 
An $n$-Lie algebra $L$ is called nilpotent of class $c$ (for some positive and integer number $c$), if $L^{c+1}=0$ and $L^c\neq0$. This is equivalent to $Z_{c-1}(L)\subsetneq Z_c(L)=L$. Then $c$ is said the nilpotency class of $L$, and we write $cl(L)=c$. Note that the nilpotency property of $n$-Lie algebras is closed under subalgebra, ideal, and homomorphism image, but it is not closed under the extension property.  

For more information about nilpotency and solvability, we refer the interested reader to  \cite{Eghdami-Gholami}. 

Let $L$ be an $n$-Lie algebra over a field $\F$ with a free presentation
\[0\lra R\lra F\lra L\lra0,\]
where $F$ is a free $n$-Lie algebra. Then the $c$-nilpotent multiplier $\M^{(c)}(L)$ of $L$ is defined as
\[\M^{(c)}(L):=\frac{R\cap F^{c+1}}{\gamma_{c+1}[R,F,\ldots,F]}.\]
So far, several studies have been done in the case $n=2$, that is, for Lie algebras. For more information, we refer to \cite{Bosko, Eshrati-Saeedi-Darabi,Niroomand-Parvizi-M^2(L), Niroomand-Russo, Salemkaretal,Salemkar-Edalatzadeh-Araskhan, c-multiplier}. 

The  authors \cite{Akbarossadat-Saeedi-3} introduced the wedge (exterior) product of two $n$-Lie algebra. Also, 
The  authors \cite{Akbarossadat-Saeedi1}  defined the non-abelian tensor products of two $n$-Lie algebras and the modular $n$-tensor products as follows.
\begin{definition}[Modular $n$-tensor product/$n$-tensor spaces]
Let $n\geq1$ be a natural number. Let $V_1$ and $V_2$ be two vector spaces over a field $\F$ of finite dimensions $d_1$ and $d_2$, respectively. Also, let $V_1^{\times i}\times V_2^{\times (n-i)}$ denote the Cartesian product
\[\underbrace{V_1\times\cdots\times V_1}_{i\ \text{times}}\times\underbrace{V_2\times\cdots\times V_2}_{n-i\ \text{times}}\]
for all $1\leq i\leq n-1$.

A function $f$ from $V_1^{\times i}\times V_2^{\times (n-i)}$ to a vector space $W$ is multilinear (or $n$-linear) if the restriction of $f$ on every component of $V_1^{\times i}\times V_2^{\times (n-i)}$ is linear. \\ 
Let $\{e_{ij}:1\leq j\leq d_i\}$ be a basis of $V_i$ for $i=1,2$. Then there exists a unique multilinear function $f:V_1^{\times i}\times V_2^{\times (n-i)}\lra W$ admitting the legal values on the elements of 
\begin{multline}\label{multilinear}
\{(e_{1 j_1},\ldots,e_{1 j_i},e_{2 j_{i+1}},\ldots,e_{2 j_n}):\\
1\leq j_k\leq d_1,\ 1\leq k\leq i,\ 1\leq j_s\leq d_2,\ i+1\leq s\leq n\}.
\end{multline}
Note that the above set contains $d_1^i\times d_2^{n-i}$ elements, which exceeds the dimension of $V_1^{\times i}\times V_2^{\times (n-i)}$ most of the time.

A pair $(\T_i,\Phi_i)$ (where $\T_i$ is a vector space and $\Phi_i$ is a multilinear function from $V_1^{\times i}\times V_2^{\times (n-i)}$ to $\T_i$) satisfies the universal factorization property if for each vector space $W$ and an $n$-linear function $f:V_1^{\times i}\times V_2^{\times (n-i)}\lra W$, there is a linear function $h_i:\T_i\lra W$ such that $f=h_i\Phi_i$. The existence of the universal pair can be proved easily. Also, up to isomorphism, there exists a unique universal pair $(\T_i,\Phi_i)$ satisfying the universal factorization property. Hence, $\T_i$ is a modular $i$-tensor product, which we may denote it by 
\[\underbrace{V_1\otimes\cdots\otimes V_1}_{i\ \text{times}}\otimes\underbrace{V_2\otimes\cdots\otimes V_2}_{n-i\ \text{times}},\]
or simply, by $V_1^{\otimes i}\otimes V_2^{\otimes (n-i)}$ for all $1\leq i\leq n-1$. Now, let
\[V\otimes_{\mathrm{mod}}^n W=\spn\{\T_i:1\leq i\leq n-1\}.\]
The vector space $V\otimes_{\mathrm{mod}}^n W$ is called the modular tensor product (or the abelian tensor product) of $V$ and $W$. 

It is evident that $\T_i=\T_j=V\otimes_{\mathrm{mod}}^n V$ for all $1\leq i,j\leq n-1$ whenever $V=W$, and that $V\otimes_{\mathrm{mod}}^n W$ coincides with the ordinary tensor product of two vector spaces when $n=2$.
\end{definition}
\begin{remark}\label{dim-modular-tensor}
For every finite-dimensional vector spaces $V$ and $W$ with $\dim V=d_V$ and $\dim W=d_W$, we have 
\[\dim(V\otimes_{\mathrm{mod}}^n W)=\sum_{i=1}^{n-1}d_V^id_W^{n-i},\]
and also, 
\[\dim((V\otimes_{\mathrm{mod}}^n V)\otimes_{\mathrm{mod}}^n W)=\sum_{i=1}^{n-1}d_V^{ni}d_W^{n-i}.\] 
\end{remark}
\bd[$c$-Capable $n$-Lie algebras]
An $n$-Lie algebra $L$ is called $c$-capable, if there exists an $n$-Lie algebra $H$ such that $L\cong H/Z_c(H)$. 
\ed 
\bd 
Let $L$ be an $n$-Lie algebra. We define $Z_c^*(L)$ to be the smallest ideal as $M$ of $L$ such that the $n$-Lie algebra $L/M$ is $c$-capable.  
\ed 
It is easy to check that $Z_c^*(L)$ is a characteristic ideal of $L$ and also, that $Z_c^*(L/Z_c^*(L))=0$.
The notation $Z^*(L)$  was defined for Lie algebras in \cite{Salemkaretal}, and it was shown that  a Lie algebra $L$ is capable if and only if  $Z^*(L)=0$. Similarly, it can be proved for $n$-Lie algebras. 
\bl \label{Z^*=0}
Let $L$  be an $n$-Lie algebra. Then $Z^*(L)=0$ if and only if $L$ is capable.
\el 
\bpr 
Let $Z^*(L)=0$. By the definition of $Z^*(L)$ we know that $L/Z^*(L)$ is capable, that is, there is an $n$-Lie algebra $M$ such that $L/Z^*(L)\cong M/Z(M)$. Since $Z^*(L)=0$, so $L/Z^*(L)\cong L$ and hence $L$ is capable. 

On the other hand, if $L$ is capable, then there exists an $n$-Lie algebra $M$ such that $L\cong M/Z(M)$. Since $Z^*(L)$ is the smallest ideal of $L$ such  that $L/Z^*(L)$ is capable,  we obtain that $Z^*(L)=0$. 
\epr 
The following theorem was proved for Lie algebras in \cite{c-multiplier}. 
\bt \label{natural epimorphism} 
Let $L$ be an $n$-Lie algebra with the free representation $F/R$. If the map $\pi:F/\gamma_{c+1}(R,F,\dots,F)\longrightarrow F/R$ is an epimorphism, then $Z_c^*(L)=\pi\left(Z_c(F/\gamma_{c+1}(R,F,\dots,F))\right)$.
\et 
\bpr 
To prove this, we have to prove that $L/\pi\left(Z_c(F/\gamma_{c+1}(R,F,\dots,F))\right)$ is $c$-capable. Also, we show that it is the smallest ideal of $L$ such that $L/\pi\left(Z_c(F/\gamma_{c+1}(R,F,\dots,F))\right)$ is $c$-capable. 
We know that $Z_c(F/\gamma_{c+1}(R,F,\dots,F))=\dfrac{Z_c(F)+\gamma_{c+1}(R,F,\dots,F)}{\gamma_{c+1}(R,F,\dots,F)}$, and hence $\pi(Z_c(F/\gamma_{c+1}(R,F,\dots,F)))=Z_c(F)/R$. Therefore, by choosing $H=F$, we have 
\[\dfrac{H}{Z_c(H)}=\dfrac{F}{Z_c(F)}\cong \dfrac{F/R}{Z_c(F)/R}=\dfrac{L}{\pi\left(Z_c(F/\gamma_{c+1}(R,F,\dots,F))\right)}.\]
Let $M$ be an ideal of $L$ such that $L/M$ is $c$-capable and $M\subseteq \pi\left(Z_c(F/\gamma_{c+1}(R,F,\dots,F))\right)$. Since $L/M$ is $c$-capable,  there is an $n$-Lie algebra $N$ such that $L/M\cong N/Z_c(N)$. We have $\pi\left(Z_c(F/\gamma_{c+1}(R,F,\dots,F))\right)/M\unlhd L/M$, since $M\subseteq \pi\left(Z_c(F/\gamma_{c+1}(R,F,\dots,F))\right)$. Therefore, 
\[\dim\dfrac{L}{\pi\left(Z_c(F/\gamma_{c+1}(R,F,\dots,F))\right)}=\dim\dfrac{L/M}{\pi\left(Z_c(F/\gamma_{c+1}(R,F,\dots,F))\right)/M}\leq \dim N/Z_c(N),\] 
and hence $M=\pi\left(Z_c(F/\gamma_{c+1}(R,F,\dots,F))\right)$. Thus $Z_c^*(L)=\pi\left(Z_c(F/\gamma_{c+1}(R,F,\dots,F))\right)$. 
\epr 
\bl \label{lem2.1}
 Let $L$ be an $n$-Lie algebra and let $N$ be an ideal of $L$. Then
$N\subseteq Z^*(L)$ if and only if the natural map $\mathcal{M}(L)\longrightarrow\mathcal{M}(L/N)$ 
is a  monomorphism.
\el 
\bpr
According to the definition of multiplier of $n$-Lie algebras, we know that 
\[\M(L)\cong\dfrac{R\cap F^2}{\gamma_2(R,F,\dots,F)},\qquad \M(L/N)\cong\dfrac{S\cap F^2}{\gamma_2(S,F,\dots,F)},\]
where $F/R$ and $S/R$ are free representations of $L$ and $N$, respectively. Thus $R\subseteq S\unlhd F$, and hence $R\cap F^2\subseteq S\cap F^2$. If $N\subseteq Z^*(L)$, then by Theorem \ref{natural epimorphism}, 
\[S/R\subseteq \pi(Z(F/\gamma_2(R,F,\dots,F))=\pi((Z(F)+\gamma_2(R,F,\dots,F))/\gamma_2(R,F,\dots,F)).\] 
Therefore, the natural map $\alpha:\mathcal{M}(L)\longrightarrow\mathcal{M}(L/N)$ 
is a  monomorphism. 

Now conversely, let $\alpha$ be a monomorphism and let  $N\not\subseteq Z^*(L)$. Then there exists $n_0\in N$ such that $n_0\notin Z^*(L)=\pi(Z(F/\gamma_2(R,F\dots,F)))$. Assume that $n_0=s_0+R\in S/R=N$, for some $s_0\in S$. Then $0\neq \pi^{-1}(s_0+R)=s_0+\gamma_2(R,F,\dots,F)\not\in Z(F/\gamma_2(R,F,\dots,F))$. Therefore, there are $f_i+\gamma_2(R,F,\dots,F)\in F/\gamma_2(R,F,\dots,F)$, for $i=2,\ldots,n$, such that $\bar{x}=[s_0,f_2,\dots,f_n]+\gamma_2(R,F,\dots,F)\neq 0$. Now, there are two cases: First, if $[s_0,f_2,\dots,f_n]\not\in R$, then $\alpha(\bar{x})$ is not defined, and hence the statement is hold. Second, if $[s_0,f_2,\dots,f_n]\in R$, then 
\[\alpha(\bar{x})=\alpha([s_0,f_2,\dots,f_n]+\gamma_2(R,F,\dots,F))=[s_0,f_2,\dots,f_n]+\gamma_2(S,F,\dots,F)=0,\]
and this contradicts with $\ker\alpha=0$.  
\epr
\bt[\cite{capability-darabi}] \label{abelian-capable}
The $d$-dimensional abelian $n$-Lie algebra $F(d)$ is capable if and only if $d\geq n$.
\et 
The following proposition is similar to Lie algebras. 
\bp \label{Z* subset I}
Let  $L$ be an $n$-Lie algebra and let $M$ be an ideal of it, such that $L/M$ is capable. Then
$Z^*(L)\subseteq M$.
\ep 
\bpr 
To prove this, it is enough to show that $Z^*(L)=\bigcap\limits_{i=1}^r M_i$, where $M_i$'s are ideals of $L$ such that $L/M_i$ is capable, for each $i=1,2,\dots,r$. It is clear that $Z^*(L)\subseteq\bigcap\limits_{i=1}^r M_i$, because $Z^*(L)$ is the smallest ideal of $L$ such that $L/Z^*(L)$ is capable. Hence $Z^*(L)\subseteq M_i$, for all $i$. 
On the other hand, since $L/Z^*(L)$ is capable, $Z^*(L)=M_j$, for some $1\leq j\leq r$. Thus $\bigcap\limits_{i=1}^r M_i\subseteq M_j=Z^*(L)$. Therefore, $Z^*(L)=\bigcap\limits_{i=1}^r M_i$.
\epr
\bp\label{Z* subset L^2}
Let $L$ be a nonabelian nilpotent  $n$-Lie algebra of finite dimension. Then
$Z^*(L)\subseteq  L^2$.
\ep
\bpr 
We know that $L/L^2$ is abelian with $\dim L/L^2\geq n$. Hence by Theorem \ref{abelian-capable}, $L/L^2$ is capable and so by Proposition \ref{Z* subset I}, we have $Z^*(L)\subseteq L^2$. 
\epr 
The proof of the following lemmas is omitted because of their similarity with Lie algebraic states.
The exterior center of $n$-Lie algebra of $L$, denoted by $Z^\wedge(L)$, is defined by 
\[Z^\wedge(L)=\{ l\in L;~l\wedge l_2\wedge\dots\wedge l_n=0,~\text{ for all } l_2,\dots,l_n\in L\}.\]
\bl \label{lem2.2}
 Let $L$ be an $n$-Lie algebra and let $N$ be a central ideal of $L$. Then
\begin{equation*}
\mathcal{M}(L)\stackrel{\alpha}{\longrightarrow} \mathcal{M}(L/N)\stackrel{\beta}{\longrightarrow}N\cap L^2\longrightarrow 0.
\end{equation*}
\el 
\bpr
The homomorphism $\alpha$ is defined in Lemma \ref{lem2.1}. Assume that $F/R$ and $S/R$ are free representations of $L$ and $N$, respectively, where $R\subseteq S\unlhd F$. Then according to the definition of $c$-nilpotent multiplier, we have $\M(L/N)=(S\cap F^2)/\gamma_2(S,F,\dots,F)$ and $N\cap L^2=S/R\cap (F^2+R)/R$. Consider 
\[\beta:\mathcal{M}(L/N)=(S\cap F^2)/\gamma_2(S,F,\dots,F)\longrightarrow N\cap L^2=S/R\cap (F^2+R)/R\] 
by $s+\gamma_2(S,F,\dots,F)\longmapsto s+R$. It is easy to check that $\beta$ is an epimorphism and also that  $\ker \beta=\mathrm{Im}\alpha$.
\epr
\bc 
 Let $L$ be an $n$-Lie algebra. Then
\begin{equation*} 
\mathcal{M}(L)\longrightarrow \mathcal{M}(L^{ab}){\longrightarrow} L^2\longrightarrow 0.
\end{equation*}
\ec 
\bpr 
It is enough to put $N=L^2$ in Lemma \ref{lem2.2}. 
\epr 
\bp 
Let $L$ be a finite-dimensional $n$-Lie algebra and let $N$ be a central ideal of $L$. Then
\begin{enumerate}
\item 
$\dim(\M(L/N))\leq \dim(\M(L))+\dim(L^2\cap N)$ and
\item 
$\dim(\M(L/N)) = \dim(\M(L))+\dim(L^2\cap N)$ 
if and only if 
$N\subseteq Z^*(L)$. 
\end{enumerate}
\ep 
\bpr 
According to Lemma \ref{lem2.2}, we have 
\begin{equation*}
\mathcal{M}(L)\stackrel{\alpha}{\longrightarrow} \mathcal{M}(L/N){\longrightarrow}N\cap L^2\longrightarrow 0.
\end{equation*}
Since $\alpha$ is not a monomorphism, so 
$\mathcal{M}(L/N)\subseteq \mathcal{M}(L)\oplus N\cap L^2$, and hence 
\[\dim(\M(L/N))\leq \dim(\M(L))+\dim(L^2\cap N).\]
If $N\subseteq Z^*(L)$, then by Lemma \ref{lem2.1}, $\alpha$ is a monomorphism. Hence 
\begin{equation*}
0\longrightarrow \mathcal{M}(L)\stackrel{\alpha}{\longrightarrow} \mathcal{M}(L/N){\longrightarrow}N\cap L^2\longrightarrow 0.
\end{equation*}
Also, $\mathcal{M}(L/N)=\mathcal{M}(L)\oplus N\cap L^2$. Therefore, 
\[\dim(\M(L/N)) = \dim(\M(L))+\dim(L^2\cap N).\]
\epr 
\bc \label{cor2.3}
Let $L$ be an $n$-Lie algebra and let $N$ a central ideal of $L$. It holds that $N \subseteq Z^\wedge(L)$ if and only if the natural map 
$L\wedge L\longrightarrow L/N \wedge L/N$ 
is a monomorphism.
\ec 
\bpr 
Employing
$L\wedge L\cong \dfrac{L\otimes L}{L\square L}$
 and 
\cite[Corollary 3.2]{Akbarossadat-Saeedi2} implies that there exists an exact series $N\otimes L\longrightarrow L\otimes L\longrightarrow L/N \otimes L/N\longrightarrow 0$ such that
$N\subseteq L^2$. 
The definition of the nonabelian exterior product completes the proof.
\epr
\bl \label{Z^=Z*}
Let $L$ be a finite-dimensional $n$-Lie algebra. Then
$Z^\wedge(L)=Z^*(L)$. 
\el 
\bpr 
It can be proved similar to \cite{Niroomand-parvizi-russo}, which was proved for Lie algebras. 
\epr 
\bc 
Let $L$ be an $n$-Lie algebra and let $N$ be an one-dimensional central ideal of $L$. Then $L$ is capable if and only if 
$\M(L)\longrightarrow \M(L/N)$ has a nontrivial kernel. 
\ec
\bpr 
It is easy to prove from Lemmas \ref{Z^*=0}, \ref{lem2.1}, and \ref{Z^=Z*}. 
\epr
The Heisenberg Lie algebras are the ubiquities in the theory of Lie algebras, in particular nilpotent Lie algebras and their classifications. An $n$-Lie algebra $L$ is called Heisenberg if $L^2=Z(L)$ and $\dim(L^2)=1$.
\begin{theorem}[\cite{Eshrati-Saeedi-Darabi}]\ 
\begin{itemize}
\item[(a).]Every Heisenberg Lie algebra is isomorphic to the following Lie algebra of odd dimension:
\[H(2,m)=\gen{x,x_1,\ldots,x_{2m}:[x_{2i-1},x_{2i}]=x,i=1,\ldots,m}.\]
\item[(b).]Let $H(n,m)$ be a Heisenberg $n$-Lie algebra of dimension $mn+1$. Then
\[\dim(\M(H(n,m)))=\begin{cases}n,&m=1,\\\\\binom{mn}{n}-1,&m>1.\end{cases}\]
\end{itemize}
\end{theorem}
The following lemma establishes the structure of every finite-dimensional nilpotent $n$-Lie algebra $L$ satisfying $\dim(L^2)=1$.
\begin{lemma}[\cite{Eshrati-Saeedi-Darabi}]\label{L=H+A}
Let $L$ be a nilpotent $n$-Lie algebra of dimension $d$ with $\dim(L^2)=1$. Then there exists $m\geq1$ such that
\[L\cong H(n,m)\oplus F(d-mn-1),\]
where $F(d-mn-1)$ is an abelian $n$-Lie algebra of dimension $d-mn-1$. 
\end{lemma}
\bt[\cite{capability-darabi}] \label{Heisenberg-capable} 
The $n$-Lie algebra $H(n,m)\oplus F(k)$ with dimension $d=mn+k+1$ is capable if and only if $m=1$.
\et 
\bl \label{Z^(L)=L^2=Z(L)} 
Let $H(n,m)$ be a Heisenberg $n$-Lie algebra with $m>1$. Then 
\[Z^\wedge(H(n,m))=H^2(n,m)=Z(H(n,m)).\] 
\el  
\bpr
Since $m>1$,  by Theorem \ref{Heisenberg-capable}, $H(n,m)$ is not capable. Hence by  Lemma \ref{Z^*=0}, $Z^*(H(n,m))\neq0$. On the other hand, by Lemmas \ref{Z* subset L^2} and \ref{Z^=Z*}, we have $Z^\wedge(H(n,m))=Z^*(H(,n,m))\subseteq H^2(n,m)=Z(H(n,m))$, because $H(n,m)$ is nilpotent with dimension $mn+1<\infty$. Also, by \cite[Proposition 4.1]{Akbarossadat-Saeedi2}, we know that 
\[H(n,m)\otimes H(n,m)\cong \dfrac{H(n,m)}{H^2(n,m)}\otimes \dfrac{H(n,m)}{H^2(n,m)},\] 
and so $H(n,m)\wedge H(n,m)\cong \dfrac{H(n,m)}{H^2(n,m)}\wedge \dfrac{H(n,m)}{H^2(n,m)}$. Thus by Corollary \ref{cor2.3}, we obtain 
\[Z^\wedge(H(n,m))=Z^*(H(n,m))=H^2(n,m)=Z(H(n,m)).\] 
\epr
For the first time, the authors \cite{Akbarossadat-Saeedi-3} introduced the concept of free $n$-Lie algebras and then in \cite{Akbarossadat-Basic} defined the concept of basic commutators of weight $w$ in $d$-dimensional $n$-Lie algebras and also proved some of  its properties and the formula to calculate the number of them. In what follows, we review it. 
\begin{theorem}[\cite{Akbarossadat-Basic}]\label{main-formula-weightw}
Let the set $X=\{x_i|x_{i+1}>x_i;~i=1,2,\dots,d\}$ be an ordered set and a basis for the free $n$-Lie algebra $F$ 
 and let $w$ be a positive integer number. Then the number of basic commutators of weight $w$ is 
\begin{equation}\label{basicweightw}
l_d^n(w)=\sum_{j=1}^{\alpha_0}\beta_{j^*}\left(\sum_{i=2}^{w-1}\alpha_i{{{{d}\choose{n-1}}}\choose{w-i}}\right),
\end{equation}
where $\alpha_0={{d-1}\choose{n-1}}$, $\alpha_i$, ($2\leq i\leq w-1$) is the coefficient of the $(i-2)$th sentence in Newton's binomial expansion $(a+b)^{w-3}$, $({\text{i.e.}}~ \alpha_i={{w-3}\choose{i-2}})$ and if ${{k-1}\choose{n-1}}+1\leq j\leq {{k}\choose{n-1}}$, (for $k=n-1,n,n+1,n+2,\dots,d-1$), then $j^*={{k-1}\choose{n-1}}+1$ and $\beta_{j^*}=(d-n-j^*+2)$.
\end{theorem}
\bc[\cite{Akbarossadat-Basic}] \label{Cor-l_n^n(w)}
If $n=d$, then 
\begin{enumerate}
\item 
$l_n^n(1)=n$. 
\item 
$l_n^n(2)={{n}\choose{n}}=1$.
\item 
$l_n^n(3)={{n}\choose{n-1}}=n$. 
\item 
$l_n^n(4)
={{n}\choose{2}}+n$.
\end{enumerate}
\ec
\bt[\cite{Akbarossadat-Basic}] \label{F^i/F^j}
Let $F$ be a free $n$-Lie algebra and let $F^i$ be the $i$th term of the lower central series of $F$, for each $i\in\Bbb N$. Then $\dfrac{F^i}{F^{i+c}}$ is the abelian of dimension $\sum\limits_{j=0}^{c-1}l_d^n(i+j)$, where $c=1,2,\ldots$.   
\et 
\section{$2$-multipliers of $n$-Lie algebras}

\indent 
Salemkar and  Aslizadeh \cite{nilpotent-multiplier of direct sum-Salemkar-Aslizadeh} presented an explicit formula for the $c$-nilpotent multipliers of the direct sum of Lie algebras whose abelianisations are of finite dimension, and under some conditions, extended it for arbitrary Lie algebras. 
They considered two Lie algebras $A$ and $B$ with dimension $d$ and $d'$, respectively. Then they made the abelian Lie algebra $\Gamma_{c+1}(A,B)$ with dimension $l_{d+d'}(c +1)-l_d(c +1)-l_{d'}(c +1)$ (the symbol $\Gamma$ does not indicate the Whitehead quadratic functor or any generalization of it). So they identified the structure of $\Gamma_{c+1}(A,B)$ and determined  its dimension. Thus they were able to calculate the $c$-nilpotent multiplier of direct sum of two Lie algebras. 

Niroomand and Parvizi \cite{2-nilpotent-Niroomand-Parvizi} chose  another method for identifying $2$-nilpotent multipliers of a direct sum of Lie algebras. First, they considered a free representation in the form $0\longrightarrow R_1\longrightarrow F_1\longrightarrow L_1\longrightarrow 0$ and $0\longrightarrow R_2\longrightarrow F_2\longrightarrow L_2\longrightarrow 0$ for each of the Lie algebras $L_1$ and $L_2$, respectively. Then, with the help of the free product of $F_1$ and $F_2$, they made a free representation $F=F_1*F_2$ for $L_1\oplus L_2$, such that $R=R_1+R_2+[F_2,F_1]$.  So, they computed the $2$-nilpotent multiplier of $L_1\oplus L_2$ in terms of $F_i$’s and $R_i$’s as follows:
\[\M^{(2)}(L_1\oplus L_2)=\dfrac{R\cap F^3}{[[R,F], F]}=\dfrac{(R_1+R_2+[F_2,F_1])\cap (F_1* F_2)^3}{[R_1+R_2+[F_2,F_1],F_1*F_2,F_1*F_2]}.\]
They showed that $\M^{(2)}(L_1)\oplus\M^{(2)}(L_2)$ is a direct summand of $\M^{(2)}(L_1\oplus L_2)$, and so $\M^{(2)}(L_1\oplus L_2)\cong\M^{(2)}(L_1)\oplus\M^{(2)}(L_2)\oplus{\bf K}$, for some subalgebra ${\bf K}$ of $\M^{(2)}(L_1\oplus L_2)$. Furthermost, they proved that $F^3=[F_2, F_1, F_1]+[F_2,F_1,F_2]+F_1^3+F_2^3$ and $\ker\alpha \equiv [F_2, F_1,F_1]+[F_2,F_1,F_2] (\mathrm{mod}[R,F,F])$, where $\alpha:\M^{(2)}(L_1\oplus L_2)\longrightarrow\M^{(2)}(L_1)\oplus\M^{(2)}(L_2)$ is an epimorphism. Also, they proved $[F_2, F_1,F_1]+[F_2,F_1,F_2]$ (${\mathrm{mod}} F^4$) is an abelian Lie algebra generated by all basic commutators of the form $[y_i,x_j,x_k]$ and $[y_r,x_s,y_t]$, where $y_i,y_r,y_t$ and $x_j,x_k,x_s$ are taken from basic sets $Y$ and $X$ of $F_2$ and $F_1$, respectively. Finally, they showed that the following equation holds:
\[{\bf K}=[F_2, F_1,F_1]+[F_2,F_1,F_2]\equiv (L_2^{ab}\otimes L_1^{ab}\otimes L_1^{ab})\oplus 
(L_2^{ab}\otimes L_1^{ab}\otimes L_2^{ab})\quad({\mathrm{mod}} [R, F, F]).\]
Indeed, here we use a new method to prove this decomposition for $n$-Lie algebras. First, we have to prove the following lemma. 
\bl \label{2-multiplier-exact sequence}
Let $L$ be an $n$-Lie algebra. Then there exists the following exact sequence:
\[0\lra \M^{(2)}(L)\lra(L\wedge L)\wedge L\lra L^3\lra 0.\]
\el 
\bpr 
Let $0\lra R\lra F\lra L\lra 0$ be a free representation of $L$ ,and so $L\cong F/R$. We know that $(L\wedge L)\wedge L\cong\left(\dfrac{F}{R}\wedge\dfrac{F}{R}\right)\wedge\dfrac{F}{R}$ and that $\M^{(2)}(L)=\dfrac{R\cap F^3}{\gamma_3[R,F,\dots,F]}$. Consider the following maps:
\begin{align*}
\begin{array}{rcl}
\Psi:\M^{(2)}(L)&\lra & \left(\dfrac{F}{R}\wedge\dfrac{F}{R}\right)\wedge\dfrac{F}{R}\\ \\
\left[\left[f_1,\dots,f_n\right],f'_2,\dots,f'_n\right]+\gamma_3\left[R,F,\dots,F\right]&\lmt & (\bar{f_1}\wedge\dots\wedge \bar{f_n})\wedge \bar{f'_2}\wedge\dots\wedge \bar{f'_n}
\end{array}
\end{align*}
\begin{align*} 
\begin{array}{rcl}
\Phi: \left(\dfrac{F}{R}\wedge\dfrac{F}{R}\right)\wedge\dfrac{F}{R}&\lra & \M^{(2)}(L)\\ \\
 (\bar{f_1}\wedge\dots\wedge \bar{f_n})\wedge \bar{f'_2}\wedge\dots\wedge \bar{f'_n}&\lmt &[[f_1,\dots,f_n],f'_2,\dots,f'_n]+\gamma_3[R,F,\dots,F] .
\end{array}
\end{align*}
It is easy to check that $\Psi$ and $\Phi$ are monomorphism and epimorphism, respectively. Also, $\Phi\Psi=0$ and $\ker\Phi\subseteq \Ima\Psi$, and hence $\ker\Phi=\Ima\Psi$. 
\epr 
The next theorem plays an essential role in proving the main results of this paper. In this theorem, we get the $2$-nilpotent multiplier of the direct sum of two $n$-Lie algebras.
\bt \label{2-multiplier of direct sum} 
Let $L$ and $M$ be two finite-dimensional $n$-Lie algebras. Then 
\begin{equation*}
\M^{(2)}(L\oplus M)\cong \M^{(2)}(L)\oplus \M^{(2)}(M)\oplus ((L^{ab}\otimes_{\mathrm{mod}}^n L^{ab}) \otimes_{\mathrm{mod}}^n M^{ab})\oplus ((M^{ab}\otimes_{\mathrm{mod}}^n M^{ab})\otimes_{\mathrm{mod}}^n L^{ab}).
\end{equation*}
\et 
\bpr 
By  Lemma \ref{2-multiplier-exact sequence}, we have 
\begin{align}
&0\lra \M^{(2)}(L\oplus M)\lra((L\oplus M)\wedge (L\oplus M))\wedge (L\oplus M)\lra (L\oplus M)^3\lra 0,\label{exact-sequ-1}\\ 
&0\lra \M^{(2)}(L)\lra(L\wedge L)\wedge L\lra L^3\lra 0,\label{exact-sequ-2}\\
&0\lra \M^{(2)}(M)\lra(M\wedge M)\wedge M\lra M^3\lra 0.\label{exact-sequ-3}
\end{align}
On the other hand, 
\begin{align}
&((L\oplus M)\wedge (L\oplus M))\wedge (L\oplus M)\nonumber\\
&\cong \left((L\wedge L)\oplus (L^{ab}\otimes_{\mathrm{mod}}^n M^{ab})\oplus (M\wedge M)\right)\wedge(L\oplus M)\nonumber\\
&\cong ((L\wedge L)\wedge (L\oplus M))\oplus \left(\left(L^{ab}\otimes_{\mathrm{mod}}^n M^{ab}\right)\wedge (L\oplus M)\right)\oplus \left((M\wedge M)\wedge (L\oplus M)\right)\nonumber\\
&\cong (((L\wedge L)\wedge L)\oplus ((L\wedge L)\wedge M)\oplus \left((( L^{ab}\otimes_{\mathrm{mod}}^n  M^{ab})\wedge L)\oplus ((L^{ab}\otimes_{\mathrm{mod}}^n M^{ab})\wedge M)\right)\nonumber\\
&\quad\oplus(((M\wedge M)\wedge L)\oplus ((M\wedge M)\wedge M))={\mathbf{K}}.\label{K}
\end{align}
Since $L$ and $M$ act trivially on $L^{ab}$ and $M^{ab}$, respectively (and conversely),  there are the following isomorphims:
\begin{align}
&(L^{ab}\otimes_{\mathrm{mod}}^n  M^{ab})\wedge L^{ab}\cong (L^{ab}\otimes_{\mathrm{mod}}^n M^{ab})\otimes_{\mathrm{mod}}^n L^{ab} \cong (L^{ab}\otimes_{\mathrm{mod}}^n L^{ab})\otimes_{\mathrm{mod}}^n M^{ab},\nonumber\\
&(L^{ab}\otimes_{\mathrm{mod}}^n  M^{ab})\wedge M^{ab}\cong (L^{ab}\otimes_{\mathrm{mod}}^n M^{ab})\otimes_{\mathrm{mod}}^n M^{ab} \cong (M^{ab}\otimes_{\mathrm{mod}}^n M^{ab})\otimes_{\mathrm{mod}}^n M^{ab}.\label{isomorohisms}
\end{align} 
We know that 
\begin{equation}\label{(L+M)^3}
(L\oplus M)^3\cong L^3\oplus M^3.
\end{equation} 
Moreover, we have the following exact sequences:
\begin{align}
&0\rightarrow (L^{ab}\otimes_{\mathrm{mod}}^n L^{ab})\otimes_{\mathrm{mod}}^n M^{ab}\rightarrow ((L\wedge L)\wedge M)\oplus (L^{ab}\otimes_{\mathrm{mod}}^n L^{ab})\otimes_{\mathrm{mod}}^n M^{ab}\rightarrow 0\rightarrow 0, \label{exact-sequ-6}\\
&0\rightarrow (M^{ab}\otimes_{\mathrm{mod}}^n M^{ab})\otimes_{\mathrm{mod}}^n L^{ab}\rightarrow ((M\wedge M)\wedge L)\oplus (M^{ab}\otimes_{\mathrm{mod}}^n M^{ab})\otimes_{\mathrm{mod}}^n L^{ab}\rightarrow 0\rightarrow 0. \label{exact-sequ-7}
\end{align}
Now, by calculating the direct sum of the corresponding sentences in sequences \eqref{exact-sequ-2}, \eqref{exact-sequ-3}, \eqref{exact-sequ-6}, and \eqref{exact-sequ-7} and applying the isomorphisms \eqref{isomorohisms} and equation \eqref{(L+M)^3}, we obtain the following exact sequence:
\begin{equation}\label{exact-sequ-8}
0\rightarrow {\mathbf{P}}\rightarrow\mathbf{K}\rightarrow L^3\oplus M^3\rightarrow 0,
\end{equation} 
where ${\mathbf{P}}=\M^{(2)}(L)\oplus \M^{(2)}(M)\oplus ( (L^{ab}\otimes_{\mathrm{mod}}^n L^{ab})\otimes_{\mathrm{mod}}^n M^{ab})\oplus ((M^{ab}\otimes_{\mathrm{mod}}^n M^{ab})\otimes_{\mathrm{mod}}^n L^{ab})$. \\
Since the sequences \eqref{exact-sequ-1} and \eqref{exact-sequ-8} are exact,  by comparing them with each other, we conclude that 
\begin{equation*}
\M^{(2)}(L\oplus M)\cong \M^{(2)}(L)\oplus \M^{(2)}(M)\oplus ((L^{ab}\otimes_{\mathrm{mod}}^n L^{ab}) \otimes_{\mathrm{mod}}^n M^{ab})\oplus ((M^{ab}\otimes_{\mathrm{mod}}^n M^{ab})\otimes_{\mathrm{mod}}^n L^{ab}).
\end{equation*}
\epr 
Note that since every one-dimensional $n$-Lie algebra $L$ is abelian and hence isomorphic to $A(1)$,  by the definition of $c$-nilpotent multiplier of $n$-Lie algebras, $\M^{(2)}(L)=0$. 

In what follows, we state a result that is also proved in Lie algebras.
\bt \label{2-multiplier-abelain}
Let $L$ be an abelian $n$-Lie algebra with finite dimension $d$. Then 
$\dim\M^{(c)}(L) = l^n_d(c + 1)$.
 In particular,
$\dim\M(L) =l_d^n(2)=\dfrac{1}{2}d(d-1)$.
\et 
\bpr 
Let $F$ be a free $n$-Lie algebra on $d$ elements. By Theorem \ref{F^i/F^j}, $F/F^2$ is an
abelian $n$-Lie algebra of dimension $d$, and so it is isomorphic to $L$ and hence $R=F^2$. Thus 
\[\M^{(c)}(L)=\dfrac{R\cap \gamma_{c+1}(F)}{\gamma_{c+1}(R,F)}=\dfrac{F^2\cap \gamma_{c+1}(F)}{\gamma_{c+1}(F^2,F)}=\dfrac{\gamma_{c+1}(F)}{\gamma_{c+2}(F)}.\]
Hence $\dim\M^{(c)}(L)=\dim \gamma_{c+1}(F)/\gamma_{c+2}(F)=l_d^n(c+1)$, 
which gives the result. 
\epr 
Eshrati, Saeedi, and  Darabi \cite{Eshrati-Saeedi-Darabi} proved that every finite-dimensional nilpotent $n$-Lie algebra can be decomposed into the direct sum of one Heisenberg $n$-Lie algebra and one abelian $n$-Lie algebra. 
So in order to get the $2$-nilpotent multiplier of each $n$-Lie algebra, we must first calculate the $2$-multiplier of each Heisenberg $n$-Lie algebra. In the following two theorems, we identify its $2$-nilpotent multiplier.
\bt \label{2-multiplier-H(n,1)}
Let $H(n,1)$ be a Heisenberg $n$-Lie algebra of dimension $n+1$. Then 
\[\M^{(2)}(H(n,1))\cong A\left(\dfrac{n^2+3n}{2}\right)
.\]
\et 
\bpr 
Since $H(n,1)$ is a nilpotent of class $2$,  its free representation is $F/F^3$ and so $R=F^3$, where $F$ is a free $n$-Lie algebra on $n$ letters. By the definition of $2$-nilpotent multiplier of $n$-Lie algebras, we have 
\begin{align*}
\M^{(2)}(H(n,1))=\dfrac{R\cap F^3}{\gamma_3(R,F,\dots,F)}=\dfrac{F^3}{[[F^3,F,\dots,F],F,\dots,F]}=\dfrac{F^3}{F^5}.
\end{align*} 
On the other hand, by Theorem \ref{F^i/F^j}, we know that 
$\dim(F^3/F^5)=l_n^n(3)+l_n^n(4)$, 
where $l_n^n(i)$ ($i=3,4$) is the number of basic commutators of weight $i$ on $n$ letters in free $n$-Lie algebras. Since $\M^{(2)}(H(n,1))$ is abelian,
\begin{align*}
\M^{(2)}(H(n,1))\cong A(l_n^n(3)+l_n^n(4))&=A\left(n+{{n}\choose{2}}+n\right)\\
&=A\left(2n+\dfrac{n(n-1)}{2}\right)\\
&=A\left(\dfrac{n^2+3n}{2}\right).
\end{align*}
\epr 
We know that a Heisenberg $n$-Lie algebra $H(n,m)$ is capable if and only if $m=1$ and that if $m\geq2$, then $H(n,m)$ is not capable. Thus the following theorem can be proved for noncapable Heisenberg $n$-Lie algebras. 
\bt \label{2-multiplier-H(n,m)}
Let $H(n,m)$ be a Heisenberg $n$-Lie algebra of dimension $mn+1$ with $m\geq2$. Then
\[\M^{(2)}(H(n,m))\cong A(l_{mn}^n(3))
.\] 
\et 
\bpr 
Since $m\geq2$, so $H(n,m)$ is not capable, and hence by Lemma \ref{Z^(L)=L^2=Z(L)}, we have $Z^\wedge(H(n,m))=H^2(n,m)=Z(H(n,m))$. 
On the other hand, if put $I=Z^\wedge(H(n,m))$, then by Lemma \ref{lem2.3}, we have 
$\M^{(2)}(H(n,m))\cong \M^{(2)}(H(n,m)/H^2(n,m))$. Since $H(n,m)/H^2(n,m)$ is abelian of dimension $mn$,  by Theorem \ref{2-multiplier-abelain} for $c=2$, we have 
\[\M^{(2)}(H(n,m)/H^2(n,m))\cong A(l_{mn}^n(3))
.\]  
Therefore, 
\[\M^{(2)}(H(n,m))\cong A(l_{mn}^n(3))
.\] 
\epr 
In the following theorem, we determine the $2$-nilpotent multiplier of $n$-Lie algebras with one-dimensional derived subalgebra.
\bt \label{2-multiplier-n-Lie algebras}
Let $L$ be an $n$-Lie algebra of dimension $d$ and with $\dim L^2=m$. Then 
\[\M^{(2)}(L)=\begin{cases}
{{A\left(\dfrac{n^2+3n}{2}+l_{d-n-1}^n(3)+\displaystyle\sum_{i=1}^{n-1}n^{ni}(d-n-1)^{n-i}+(d-n-1)^{ni}n^{n-i}\right)}}, \\ 
\qquad\hspace{11.5cm} m=1, \\ 
{{A\left(l_{mn}^n(3)+l_{d-mn-1}^n(3)+\displaystyle\sum_{i=1}^{n-1}(mn)^{ni}(d-mn-1)^{n-i}+(d-mn-1)^{ni}(mn)^{n-i}\right)}} ,\\
\qquad\hspace{11.5cm} m\geq2.
\end{cases}\]
\et 
\bpr 
Since $\dim L^2=1$, by Lemma \ref{L=H+A}, $L\cong H(n,m)\oplus A(d-mn-1)$. Thus by Lemma \ref{2-multiplier of direct sum}, we have 
\begin{align}
\M^{(2)}(L)&\cong \M^{(2)}(H(n,m)\oplus A(d-mn-1))\nonumber\\
&\cong \M^{(2)}(H(n,m))\oplus \M^{(2)}(A(d-mn-1))\nonumber\\
&\quad\oplus \left((H^{ab}(n,m)\otimes_{\mathrm{mod}^n}A(d-mn-1))\otimes_{\mathrm{mod}^n}H^{ab}(n,m)  \right)\nonumber\\ 
&\quad\oplus \left((H^{ab}(n,m)\otimes_{\mathrm{mod}^n} A(d-mn-1))\otimes_{\mathrm{mod}^n} A(d-mn-1)\right).\label{equ1}
\end{align}
Now, we consider two following cases: 
\begin{enumerate}
\item[(a).] 
Assume that $m=1$. By Theorems \ref{2-multiplier-abelain} and \ref{2-multiplier-H(n,1)}, we have 
\begin{align*}
\M^{(2)}(L)&\cong A(l_n^n(3)+l_n^n(4)) \oplus A(l_{d-n-1}^n(3))\\
&\quad\oplus A\left(\sum_{i=1}^{n-1}n^{ni}(d-n-1)^{n-i}\right)\\
&\quad\oplus A\left(\sum_{i=1}^{n-1}(d-n-1)^{ni}n^{n-i}\right) \\
&\cong A\left(l_n^n(3)+l_n^n(4)+l_{d-n-1}^n(3)+\sum_{i=1}^{n-1}n^{ni}(d-n-1)^{n-i}+(d-n-1)^{ni}n^{n-i}\right)\\
& \cong A\left(\dfrac{n^2+3n}{2}+l_{d-n-1}^n(3)+\sum_{i=1}^{n-1}n^{ni}(d-n-1)^{n-i}+(d-n-1)^{ni}n^{n-i}\right).
\end{align*}
\item[(b).]
Now, suppose that $m\geq2$. By Theorems \ref{2-multiplier-abelain} and \ref{2-multiplier-H(n,m)}, we have 
\begin{align*}
\M^{(2)}(L)&\cong A\left(l_{mn}^n(3)
\right) \oplus A\left(l_{d-mn-1}^n(3)\right)\\
&\quad\oplus A\left(\sum_{i=1}^{n-1}(mn)^{ni}(d-mn-1)^{n-i}\right)\\
&\quad\oplus A\left(\sum_{i=1}^{n-1}(d-mn-1)^{ni}(mn)^{n-i}\right) \\
&\cong A\left(l_{mn}^n(3)+l_{d-mn-1}^n(3)+\sum_{i=1}^{n-1}(mn)^{ni}(d-mn-1)^{n-i}+(d-mn-1)^{ni}(mn)^{n-i}\right).
\end{align*}
\end{enumerate}
\epr 
\bp \label{Q_3}
Let $L$ be an $n$-Lie algebra with free representation $F/R$ and let $M$ be an ideal of $L$, such that $M\cong S/R$, for some ideal $S$ of $F$. Then $Q_3=\dfrac{\gamma_3(S,F,\dots,F)}{\gamma_3(R,F,\dots,F)}$ is an isomorphic image of $(M\wedge L)\wedge L$. Moreover, if $M$ is an $r$-central ideal, that is, $M\subseteq Z_r(L)$, then $Q_3$ is the  isomorphic image of $(M\wedge K)\wedge K$, where $K=L/\gamma_{r+1}(L)$. 
\ep 
\bpr 
We know that 
\[Q_3=\dfrac{\gamma_3(S,F,\dots,F)}{\gamma_3(R,F,\dots,F)}=\dfrac{[[
S,F,\dots,F],F,\dots,F]}{[[R,F,\dots,F],F,\dots,F]}.\]
Since $L\cong F/R$, there exists an epimorphism $g:F\longrightarrow L$, with $\ker g=R$, and hence $\bar{g}:F/R\longrightarrow L$ is an isomorphism. Consider then map $\alpha:(M\wedge L)\wedge L\rightarrow Q_3$, such that $\alpha$ maps every element $\left(\sum_{j=1}^{n-1}(m_1\wedge \dots\wedge m_j\wedge l_{j+1}\wedge \dots\wedge l_n)\wedge l_{i+1}'\wedge \dots\wedge l_n'\right)$ from $(M\wedge L)\wedge L$ to $\sum_{j=1}^{n-1}[[s_1,\dots,s_j,f_{j+1},\dots,f_n],f_{i+1}',\dots,f_n']+\gamma_3(R,F,\dots,F)$, 
where $\bar{g}^{-1}(m_i)=s_i+R$, $\bar{g}^{-1}(l_{k})=f_k+R$, for $1\leq i\leq j$, $i+1\leq k\leq n$, $1\leq j\leq n-1$. Since $\bar{g}$ is an isomorphism and the bracket of $n$-Lie algebras is an $n$-linear, $\alpha$ is a well-defined homomorphism of $n$-Lie algebras. Also, it is easy to check that $\alpha$ is onto.  \\
The second part can be proved with a similar argument.
\epr 
The following proposition is an important tool to prove the main theorem of this section.
\bp \label{lem2.3}
Let $L$ be an $n$-Lie algebra and let $M$ be an ideal of it. Then 
\begin{equation}\label{inequality-(1)}
\dim\M^{(2)}(L/M)\leq \dim\M^{(2)}(L)+\dim\dfrac{M\cap L^2}{\gamma_3(M,L,\dots, L)}.
\end{equation}
Moreover, if $M$ is a $2$-central subalgebra (i.e., $M\subseteq Z_2(L)$), then 
\begin{align}
&(M\wedge L)\wedge L\longrightarrow\M^{(2)}(L)\longrightarrow \M^{(2)}(L/M)\longrightarrow M\cap L^3\longrightarrow 0,\label{inequality-(2)}\\
&\dim\M^{(2)}(L)+\dim M\cap L^3\leq \dim\M^{(2)}(L/M)+\dim (M\otimes L/L^3)\otimes L/L^3.\label{inequality-(3)}
\end{align}
\ep 
\bpr 
By attention to the notations of Proposition \ref{Q_3} and $Q_3\subseteq \M^{(2)}(L)$, we have the following exact sequence:
\begin{equation*}
(M\wedge L)\wedge L\stackrel{\alpha}{\longrightarrow} \M^{(2)}(L)\stackrel{\beta}{\longrightarrow} \M^{(2)}(L/M)\stackrel{\gamma}{\longrightarrow} \dfrac{M\cap L^3}{\gamma_3(M,L,\dots,L)}\longrightarrow 0,
\end{equation*} 
where $\alpha$ is defined in Proposition \ref{Q_3}. Since $R\subseteq S$, so 
\[R\cap F^3\subseteq S\cap F^3,\qquad \gamma_3(R,F,\dots,F)\subseteq \gamma_3(S,F,\dots,F).\] 
Hence $\beta:\M^{(2)}(L)=\dfrac{R\cap F^3}{\gamma_3(R,F,\dots,F)}\longrightarrow \M^{(2)}(L/M)=\dfrac{S\cap F^3}{\gamma_3(S,F,\dots,F)}$. It is easy to check that $\beta$ is a Lie homomorphism with $\ker\beta=Q_3$. Also, since $\M^{(2)}(L/M)=S\cap F^3/\gamma_3(S,F,\ldots,F)$, we define $\gamma:\M^{(2)}(L/M)\longrightarrow M\cap L^3/\gamma_3(M,L,\ldots,L)$ such that $s+\gamma_3(S,F,\dots,F)\longmapsto g(s)+\gamma_3(M,L,\dots,L)$. It is obvious that since $g$ is an epimorphism so is $\gamma$ as well. Therefore, 
the inequality \eqref{inequality-(1)} is obtained. Furthermore, if $M$ is $2$-central, then $\gamma_3(M,L,\dots,L)=0$ and $(M\wedge L)\wedge L\cong (M\otimes L/L^3)\otimes L/L^3$. Hence we have the exact sequence \eqref{inequality-(2)} and the inequality \eqref{inequality-(3)}.   
\epr 
The next theorem is the main result of this section. 
\bt \label{2-multiplier-dimL^2=k>1}
Let $L$ be a nilpotent $n$-Lie algebra of dimension $d$, with $\dim L^2=k\geq 1$. Then 
\[\dim\M^{(2)}(L)\leq \dfrac{n^2+3n}{2}+l_{d-n-k}^n(3)+\displaystyle\sum_{i=1}^{n-1}n^{ni}(d-n-k)^{n-i}+(d-n-k)^{ni}n^{n-i}+(d-k)^{2n-2}-k+1.\] 
\et 
\bpr 
In Theorem \ref{2-multiplier-n-Lie algebras}, it can be seen that in the first case, that is, $m=1$, the dimension of $\M^{(2)}(L)$ is larger than the second case, that is, $m>1$. Also, in the case $m>1$, the dimension of $\M^{(2)}(L)$ is decreasing with respect to $m$. Now, we do induction on $k$. For $k=1$, the result follows from Theorem \ref{2-multiplier-n-Lie algebras}. Now, let $k=2$ and let $I$ be a one-dimensional central ideal of $L$. Then $\dim L/M=d-1$ and $\dim (L/M)^2=\dim L^2/M=1$. Also, since $M\subseteq Z(L)\subseteq Z_2(L)$, thus $M$ and $L/L^3$ act on each other trivially. Hence $(M\otimes L/L^3)\otimes L/L^3\cong (M\otimes_{\mathrm{mod}}^n (L/L^3)^{ab})\otimes_{\mathrm{mod}}^n (L/L^3)^{ab}$. So by equation \eqref{inequality-(3)} of Proposition \ref{lem2.3} and the first step of induction, we have 
\begin{align*}
\dim\M^{(2)}(L)&\leq \dim\M^{(2)}(L/M)+\dim (M\otimes L/L^3)\otimes L/L^3-\dim M\cap L^3\\
&\leq \dim\M^{(2)}(L/M)+\dim (M\otimes_{\mathrm{mod}}^n (L/L^3)^{ab})\otimes_{\mathrm{mod}}^n (L/L^3)^{ab}-\dim M\cap L^3\\
&\leq \dim\M^{(2)}(L/M)+\dim (M\otimes_{\mathrm{mod}}^n \dfrac{L/L^3}{(L/L^3)^2})\otimes_{\mathrm{mod}}^n \dfrac{L/L^3}{(L/L^3)^2}-\dim M\cap L^3\\
&\leq \dim\M^{(2)}(L/M)+\dim (M\otimes_{\mathrm{mod}}^n \dfrac{L/L^3}{L^2/L^3})\otimes_{\mathrm{mod}}^n \dfrac{L/L^3}{L^2/L^3}-\dim M\cap L^3\\
&\leq \dim\M^{(2)}(L/M)+\dim (M\otimes_{\mathrm{mod}}^n L/L^2)\otimes_{\mathrm{mod}}^n L/L^2-\dim M\cap L^3\\
&\leq \dfrac{n^2+3n}{2}+l_{(d-1)-n-1}^n(3)\\
&\quad+\displaystyle\sum_{i=1}^{n-1}n^{ni}((d-1)-n-1)^{n-i}+((d-1)-n-1)^{ni}n^{n-i}
+(d-2)^{2n-2}-1.
\end{align*}
Now suppose that the inequality holds for $k-1$, and let $M$ be a one-dimensional central ideal of $L$, with $\dim L/M=d-1$ and $\dim L^2=k$. Then 
\begin{align*}
\dim\M^{(2)}(L)&\leq \dim\M^{(2)}(L/M)-\dim(M\cap L^3)+\dim (M\otimes_{\mathrm{mod}}^n L/L^2)\otimes_{\mathrm{mod}}^n L/L^2\\
&\leq \dfrac{n^2+3n}{2}+l_{(d-1)-n-(k-1)}^n(3)\\
&\quad+\displaystyle\sum_{i=1}^{n-1}n^{ni}((d-1)-n-(k-1))^{n-i}+((d-1)-n-(k-1))^{ni}n^{n-i}\\
&\quad-\dim(M\cap L^3)+\dim (M\otimes_{\mathrm{mod}}^n L/L^2)\otimes_{\mathrm{mod}}^n L/L^2
\end{align*}
\begin{align*}
\qquad\qquad\qquad &\leq \dfrac{n^2+3n}{2}+l_{d-n-k}^n(3)\\
&\quad+\displaystyle\sum_{i=1}^{n-1}n^{ni}(d-n-k)^{n-i}+(d-n-k)^{ni}n^{n-i}-k+1+(d-k)^{2n-2}.
\end{align*}
\epr 
\section{$2$-Capability of $n$-Lie algebras}
\indent 
Detection of $2$-capable $n$-Lie algebras is one of the applications of $2$-nilpotent multipliers. Hence in this section, we are going to state the conditions for detecting $2$-capable $n$-Lie algebras. 
\bp \label{Pro3.1}
Let $L$ be an $n$-Lie algebra. Then $L$ is $2$-capable if and only if $Z_2^*(L)=0$.  
\ep 
\bpr 
If $L$ is $2$-capable, then there exists an $n$-Lie algebra $H$ such that $L\cong H/Z_2(H)$. Assume that $F'/R'$ is a free representation of $H$. Then $Z_2(H)\cong S/R'$, where $R'\subseteq S\unlhd F'$. Thus 
\[L\cong \dfrac{H}{Z_2(H)}\cong \dfrac{F'/R'}{S'/R'}\cong \dfrac{F'}{S'}.\]
On the other hand, since $S'\unlhd F'$, so $\gamma_3(S',F',\dots,F')\subseteq S'$, $\dfrac{F'/\gamma_3(S',F',\dots,F')}{S'/\gamma_3(S',F',\dots,F')}\cong F'/S'$, and also the map $\pi:F'/\gamma_3(S',F',\dots,F')\rightarrow F'/S'$  is the natural epimorphism. Hence by  Theorem \ref{natural epimorphism}, $Z_2^*(L)=\pi\left(Z_2(F'/\gamma_3(S',F',\dots,F'))\right)$. Therefore, it is enough to show that $Z_2(F'/\gamma_3(S',F',\dots,F')\subseteq \ker\pi$. Since $S/\gamma_3(S',F',\dots,F')=Z_2(F'/\gamma_3(S',F',\dots,F'))$, so $S'/\gamma_3(S',F',\dots,F')\subseteq \ker\pi$, and hence 
$\pi(S'/\gamma_3(S',F',\dots,F'))=0$. Thus $Z_2^*(L)=0$.  

Now, let $Z_2^*(L)=0$, and suppose that $F/R$ is a free representation of $L$. We have $\gamma_3(R,F,\dots,F)\subseteq R\unlhd F$, and hence $\pi:F/\gamma_3(R,F,\dots,F)\rightarrow F/R$ is an epimorphism with $\ker\pi=S/\gamma_3(R,F,\dots,F)$. Also, 
\begin{equation}\label{equ3}
\dfrac{F/\gamma_3(R,F,\dots,F)}{\ker \pi}\cong\dfrac{F/\gamma_3(R,F,\dots,F)}{R/\gamma_3(R,F,\dots,F)}\cong \dfrac{F}{R}=\Ima \pi=L.
\end{equation} 
Put $H=F/\gamma_3(R,F,\dots,F)$. Using the definition of $Z_2(H)$, it is easy to check that 
\begin{equation}\label{equ4}
Z_2(H)=Z_2(F/\gamma_3(R,F,\dots,F))=\dfrac{R}{\gamma_3(R,F,\dots,F)}.
\end{equation} 
Therefore, $L\cong H/Z_2(H)$ and hence is $2$-capable. 
\epr 
\bt \label{Theo3.2}
Let $L$ be an $n$-Lie algebra and let $M$ be  its ideal such that $M\subseteq Z_2^*(L)$. Then the natural homomorphism $\psi:\M^{(2)}(L)\longrightarrow \M^{(2)}(L/M)$ is one-to-one.  
\et 
\bpr 
Suppose that $F/R$ and $S/R$ are free representations of $L$ and $M$, respectively, for some $S\unlhd F$ with $R\subseteq S$. Then $L/M\cong F/S$, and we have 
\[\M^{(2)}(L)=\dfrac{R\cap F^3}{\gamma_3(R,F,\dots,F)},\qquad \M^{(2)}(L/M)=\dfrac{S\cap F^3}{\gamma_3(S,F,\dots,F)}.\]
Since $R\subseteq S$, so $R\cap F^3\subseteq S\cap F^3$ and 
\begin{equation}\label{equ-1}
\gamma_3(R,F,\dots,F)\subseteq \gamma_3(S,F,\dots,F)\subseteq S.
\end{equation}  
Hence we define the homomorphism $\psi:(R\cap F^3)/\gamma_3(R,F,\dots,F)\longrightarrow (S\cap F^3)/\gamma_3(S,F,\dots,F)$ by $\psi(x+\gamma_3(R,F,\dots,F))=x+\gamma_3(S,F,\dots,F)$. 

On the other hand, since $I\subseteq Z_2^*(L)$, so $S/R\subseteq Z_2(F/R)$, or, equivalently,
\begin{equation}\label{equ-3}
[[S,F,\dots,F],F,\dots,F]=\gamma_3(S,F,\dots,F)\subseteq R.
\end{equation} 
By equations \eqref{equ-1} and \eqref{equ-3}, we obtain 
\[\gamma_3(R,F,\dots,F)\subseteq \gamma_3(S,F,\dots,F)\subseteq R\cap S=R.\]
Therefore one of the two situations $\gamma_3(S,F,\dots,F)=R$ or $\gamma_3(R,F,\dots,F)=\gamma_3(S,F,\dots,F)$ can occur. If $\gamma_3(S,F,\dots,F)=R$, then $S=R$ must be, which means $I=S/R=0$, a contradiction. Thus $\gamma_3(R,F,\dots,F)=\gamma_3(S,F,\dots,F)$. Hence for arbitrary $x+\gamma_3(R,F,\dots,F)\in\ker\psi$, we have 
\[\psi(\gamma_3(R,F,\dots,F))=x+\gamma_3(S,F,\dots,F)=\gamma_3(S,F,\dots,F).\]
So $x\in\gamma_3(S,F,\dots,F)=\gamma_3(R,F,\dots,F)$. Thus $\ker\psi=0$. 
\epr 
In the following theorem, we prove that the necessary condition for $H(n,m)$ to be a $2$-capable $n$-Lie algebra is that $m=1$.
\bt \label{necessary condition for capabilitty}
If the Heisenberg $n$-Lie algebra $H(n,m)$ is $2$-capable, then $m=1$. 
\et 
\bpr 
Suppose that $H(n,m)$ is a $2$-capable Heisenberg $n$-Lie algebra. Then there is an $n$-Lie algebra $K$ such that $H(n,m)\cong K/Z_2(K)$, and we have 
\[H(n,m)\cong K/Z_2(K)\cong \dfrac{K/Z(K)}{Z_2(K)/Z(K)}=\dfrac{K/Z(K)}{Z(K/Z(K))}.\] 
Thus putting $P=K/Z(K)$, we conclude that $H(n,m)=P/Z(P)$. Hence $H(n,m)$ is capable, and so by Theorem \ref{Heisenberg-capable}, $m=1$. 
\epr 
In the following result, we show that there is no $2$-capable Heisenberg $n$-Lie algebras, for $n>2$. 
\bt \label{non-capability of H(n,1)}
The Heisenberg $n$-Lie algebra $H(n,1)$ is not $2$-capable, for $n\geq3$. 
\et 
\bpr 
Assume that the Heisenberg $n$-Lie algebra $H(n,1)$ is $2$-capable, and $n\neq2$. So there is an $n$-Lie algebra $H$, such that $H(n,1)\cong H/Z_2(H)$ and $Z_2^*(H(n,1))=0$. By Theorem \ref{natural epimorphism}, $\pi(Z_2(F/\gamma_3(R,F,\dots,F))=0$. That is, 
\[Z_2(F/\gamma_3(R,F,\dots,F)\subseteq \ker\pi=K/\gamma_3(R,F,\dots,F),\] 
for some ideal $K$ of $F$. Hence $Z_2(F)+\gamma_3(R,F,\dots,F)/\gamma_3(R,F,\dots,F)\subseteq K/\gamma_3(R,F,\dots,F)$, and so $Z_2(F)\subseteq K$.   

On the other hand, $Z_2^*(H(n,1))=\pi(Z_2(F/\gamma_3(R,F,\dots,F))=0$ implies that 
\[\pi(Z_2(F/\gamma_3(R,F,\dots,F))=\pi(Z_2(F)+\gamma_3(R,F,\dots,F))+\gamma_3(R,F,\dots,F)=R.\]
Since $c$-central ideals are characteristic, thus $\pi(Z_2(F))=Z_2(F)\subseteq R$. 

Therefore, $Z_2(F)\subseteq K\cap R$ and that is $Z_2(F/R)=0$, which is a contradiction.  
\epr 
According to the above discussion, we have the following corollary. 
\bc 
The Heisenberg $n$-Lie algebra $H(1)$ is the only $2$-capable Heisenberg $n$-Lie algebra.
\ec
\bpr 
It is proved that $H(2,1)=H(1)$ is $2$-capable; see \cite[Theorem 3.3]{Niroomand-Parvizi-M^2(L)}. Also, by attention to Theorems \ref{non-capability of H(n,1)} and \ref{necessary condition for capabilitty}, the proof is completed.    
\epr 
The following theorem is similar to \cite[Lemma 4.2]{Moneyhun}. Moneyhun proved that if $\dim L/Z(L)=d$, then $\dim L^2=\dim\gamma_2(L)\leq\dfrac{1}{2}d(d-1)=l_d(2)$.  
\bt 
Let $L$ be an $n$-Lie algebra such that $\dim L/Z_2(L)=d$. Then the dimension of $\gamma_{3}(L)$ is at most $l_d^n(3)$. 
\et 
\bpr 
Suppose that $F/R$ is the free representation of $L$ and that $\{x_1,x_2,\dots,x_d\}$ is a basis for $L/Z_2(L)$. Thus every element of $L/Z_2(L)$ can be represented as $z+\sum_{i=1}^dc_ix_i$, where $c_i$s are scalers. Hence 
\begin{align*}
\dim\gamma_3(L)=\dim\gamma_3(F/R)&=\dim \dfrac{F^3+R}{R}
=\dim \dfrac{F^3}{F^3\cap R}
=\dim\dfrac{F^3/F^4}{(F^3\cap R)/F^4}
\leq \dim\dfrac{F^3}{F^4}.
\end{align*}
By Theorem \ref{F^i/F^j}, we know that the basis of $F^3/F^4$ is the set of all basic commutators of weight $3$, and so $\dim\gamma_3(L)\leq l_d^n(3)$.
\epr 

\ \\
{\small 
Farshid Saeedi$^*$: Corresponding author \\ 
Department of Mathematics, Mashhad Branch, Islamic Azad University, Mashhad, Iran.\\
E-mail address: saeedi@mshdiau.ac.ir \\ \ \\
Seyedeh Nafiseh Akbarossadat \\
Department of Mathematics, Mashhad Branch, Islamic Azad University, Mashhad, Iran.\\
E-mail address: n.akbarossadat@gmail.com  \\ \ \\ 
}

\end{document}